\documentclass[twocolumn]{autart}    

\usepackage{graphicx}          
\usepackage{amsmath,amsfonts}

\newtheorem{problem}{Problem}
\newtheorem{remark}{Remark}
\newtheorem{definition}{Definition}
\newtheorem{theorem}{Theorem}

\begin{document}

\begin{frontmatter}

\title{A Systematic LMI Approach to Design Multivariable Sliding Mode Controllers\thanksref{footnoteinfo}} 

\thanks[footnoteinfo]{This paper was not presented at any IFAC 
meeting. Corresponding author P.~H.~.S.~Coutinho.}

\author[UERJ]{Pedro H. S. Coutinho}\ead{phcoutinho@eng.uerj.br},    
\author[UFAM]{Iury Bessa}\ead{iurybessa@ufam.edu.br},               
\author[UERJ]{Victor H. P. Rodrigues}\ead{ victor.rodrigues@uerj.br},  
\author[UERJ]{Tiago Roux Oliveira}\ead{tiagoroux@uerj.br}  

\address[UERJ]{Department of Electronics and Telecommunication Engineering, State University of Rio de Janeiro (UERJ), RJ 20550-900, Brazil}  
\address[UFAM]{Department of Electricity, Federal University of Amazonas (UFAM), AM 69080-900, Brazil}             


\begin{abstract}                          
This paper deals with sliding mode control for multivariable polytopic uncertain systems. We provide systematic procedures to design variable structure controllers (VSCs) and unit-vector controllers (UVCs). Based on suitable representations for the closed-loop system, we derive sufficient conditions in the form of linear matrix inequalities (LMIs) to design the robust sliding mode controllers such that the origin of the closed-loop system is globally stable in finite time. Moreover, by noticing that the reaching time depends on the initial condition and the decay rate, we provide convex optimization problems to design robust controllers by considering the minimization of the reaching time associated with a given set of initial conditions. Two examples illustrate the effectiveness of the proposed approaches.    
\end{abstract}

\end{frontmatter}

\section{Introduction}

Sliding mode control has attracted attention because of its robustness, easy implementation, and the possibility to achieve fast responses~\cite{oliveira2023sliding}. For this reason, sliding mode control approaches have been applied to several control problems, such as control of time-delay~\cite{chen2024sliding}, perturbed PDE~\cite{zhang2023finite}, switched~\cite{lian2020event}, uncertain~\cite{moreno2022multivariable}, and networked control systems~\cite{cucuzzella2018practical,rinaldi2021adaptive}.

However, tuning sliding-mode controllers can be intricate and strongly dependent on a priori knowledge of the disturbances or their upper bounds~\cite{Roy2020,ovalle2025analysis}. Indeed, the usual stability criteria may not apply to relay systems, which poses an additional difficulty in designing sliding-mode controllers. For example,  systems of the form $\dot{\sigma}(t) = A \mathrm{sgn}(\sigma(t))$ may be unstable even for a Hurwitz matrix $A$~\cite{utkin2020road}. In this regard, reference~\cite{hsu2000matrix} provided sufficient stability analysis conditions for the multivariable case using a Persidskii diagonal-type Lyapunov function, and necessary conditions could be presented only for second-order systems. 

To deal with the synthesis of sliding-mode controllers, using particular classes of Lyapunov functions allows us to obtain constructive design conditions based on semidefinite programming optimization under linear matrix inequalities (LMIs) constraints. For instance, this kind of approach eases the extension of super-twisting algorithms to multivariable systems~\cite{geromel2024multivariable,geromel2024lmi}. Moreover, LMI-based conditions have been successfully applied to design sliding-mode controllers for systems with disturbances~\cite{Roy2020}, time delays~\cite{Chuei2023}, and uncertainty~\cite{geromel2024multivariable}. In~\cite{Mei2022,Mei2022b}, the authors present LMI-based sufficient conditions for input-to-state stability and synchronization of Persidskii systems. For linear and nonlinear quasi-Lipschtiz systems, the authors in~\cite{Polyakov2011} present 
an LMI-based invariant ellipsoid method to select sliding surfaces, ensuring stability and disturbance attenuation. However, the synthesis problem for the class of systems of the form $\dot{\sigma}(t)\!=\!A \mathrm{sgn}(\sigma(t))$ is still an open problem, especially when uncertainties are concerned. This is one of the motivations for this work.

\textbf{Contributions}:
In this paper, we address the sliding mode control design problem for the class of polytopic 
uncertain systems using a traditional variable structure control (VSC) with a relay-type function and a unit vector control (UVC)~\cite{hsu2002multivariable}. 
In both cases, we propose suitable state transformations from which the closed-loop dynamics 
can be rewritten in an appropriate form that enables deriving LMI-based design conditions. In the VSC case, the conditions are based on a Persidskii diagonal-type Lyapunov function candidate. 
The robust control gains designed with the derived conditions ensure the
closed-loop equilibrium is robustly finite-time stable. 
Since the reaching time estimate depends on the initial condition and the decay rate enforced in the control design~\cite{hsu2000matrix}, 
we propose optimization problems 
to perform the optimal control design by incorporating objective functions 
related to the reaching time minimization and the 
maximization of the set of initial conditions associated with the guaranteed reaching time.

\textbf{Organization:}
The paper is organized as follows.
The results for the VSC are described in Section~\ref{sec:VSC}, where
we present the class of uncertain systems studied in this paper, the robust VSC design condition, and the optimization problem.
Section~\ref{sec:UVC} presents the proposed
robust UVC design condition and the associated optimization problem.
Two examples are presented in Section~\ref{sec:results}, where
we illustrate the effectiveness of the proposed approach
in designing robust sliding mode controllers with a guaranteed
minimized reaching time with the estimated set of initial conditions
for which such a reaching time is ensured.
Conclusions are depicted in Section~\ref{sec:conclusion}.

\textbf{Notation:}
$\mathbb{N}$ is the set of natural numbers; $\mathbb{N}_{\leq m}=\{1,\ldots,m\}$ for some $m\in\mathbb{N}$. $\mathbb{R}^n$ is the $n$-dimensional Euclidean space and $\mathbb{R}^{n \times m}$ is the set of real matrices of order $n \times m$; $\mathrm{diag}(A, B)$ is a block diagonal matrix whose elements of the main-diagonal are the matrices $A$ and $B$. 
For a given
symmetric matrix $P$, $P > 0~(<0)$ indicates that $P$ is a positive (negative) definite matrix.
For a given matrix $A$, $A^\top$ denotes its transpose, and $\lambda_{\min}(A)$ is the
smallest eigenvalue of $A$.

\section{Robust Variable Structure Control Design}
\label{sec:VSC}

Consider the following uncertain MIMO system
\begin{align}\label{eq:plant}
    \dot{\sigma}  = Bu,
\end{align}
where $\sigma \in \mathbb{R}^n$ is the state vector defined as $\sigma = (\sigma_1, \ldots, \sigma_n)$,
$u \in \mathbb{R}^m$ is the input vector.
Moreover, $B$ is constant but unknown, taking values in the set $\mathcal{B} = \mathrm{co}\{B_i\}_{i \in \mathbb{N}_{\leq N}}$.
This implies that it is possible to write 
\begin{align}
    B = \sum_{i = 1}^N \alpha_i B_i,
\end{align} 
where $\alpha = (\alpha_1,\ldots,\alpha_N)$
is the vector of uncertain parameters lying in the unit simplex 
\begin{align}
    \Lambda = \left\lbrace\alpha \in \mathbb{R}^N: \sum_{i=1}^N \alpha_i = 1, \; \alpha_i \geq 0\right\rbrace.
\end{align}
In this work, we consider the following variable structure control law
\begin{align}\label{eq:controller}
    u = K \mathrm{sgn}(\sigma),
\end{align}
where the sign function is understood in the element-wise sense, that is, $\mathrm{sgn}(\sigma) = (\mathrm{sgn}(\sigma_1), \ldots, \mathrm{sgn}(\sigma_n))$. Then, the closed-loop
system is given by
\begin{align}\label{eq:closed-loop}
    \dot{\sigma} = B K \mathrm{sgn}(\sigma).
\end{align}
The problem addressed in this section is stated as follows.
\begin{problem}\label{problem:1}
    Consider the uncertain system~\eqref{eq:plant} and the VSC law~\eqref{eq:controller}.
    Design a robust control gain $K$ such that the origin of the closed-loop system~\eqref{eq:closed-loop} is globally finite-time stable. 
\end{problem}

\begin{remark}
    A similar problem of studying the global stability of a dynamical system of the form $\dot{\sigma}(t) = A \mathrm{sgn}(\sigma(t))$, where $A \in \mathbb{R}^{n \times n}$ is a constant matrix, is recognized as a
    challenging problem, especially when it comes to deriving necessary
    and sufficient conditions. A sufficient and necessary condition was established by~\cite{hsu2000matrix} for the case of $n=2$ and $A$ precisely known. 
    The authors in~\cite{hsu2000matrix} also provided a sufficient condition
    for arbitrary $n$ using a diagonal-type Lyapunov function candidate
    inspired by the analysis of Persidskii-type systems. In this work, we tackle the robust control synthesis problem for the case of an \textit{uncertain} system as~\eqref{eq:plant}. Therefore, we aim to design the control gain $K$ of the controller~\eqref{eq:controller} such that the origin of the closed-loop system~\eqref{eq:closed-loop} is globally finite-time stable,
    as stated in Problem~\ref{problem:1}.
\end{remark}



\subsection{The proposed robust VSC design condition}

This section aims to derive LMI-based conditions for designing 
robust VSCs. For that purpose, consider
the state vector $x = (\sqrt{|\sigma_1|}, \ldots, \sqrt{|\sigma_n|})$.
In $x$-coordinates, the closed-loop system~\eqref{eq:closed-loop} can be
equivalently written as follows:
\begin{align}
    \dot{x} = \frac{1}{2} L(\sigma) B K L(\sigma) x,
\end{align}
where $L(\sigma) = \mathrm{diag}(\sigma_1/|\sigma_1|^{3/2}, \ldots, \sigma_n/|\sigma_n|^{3/2})$. This new set of coordinates allows
to write the following diagonal-type Lyapunov function candidate, inspired from~\cite{hsu2000matrix}:
\begin{align}\label{eq:proof-thm1-6}
    V(\sigma) = \sum_{i=1}^{n} p_i |\sigma_i|,
\end{align}
where $p_i > 0$, in the following quadratic form:
\begin{align}
    V(x) = x^\top P x,
\end{align}
with $P = \mathrm{diag}(p_1,\ldots,p_n)$.

Before stating the main result of this work, we consider the following definition of the equivalent control input.

\begin{definition}[see~\cite{hsu2000matrix}]\label{def:extended_u}
    Consider a Filippov solution $\sigma(t)$ of
system~\eqref{eq:closed-loop}, defined for $t$ in a certain interval $T$. Then, the extended equivalent control $u_\mathrm{eq}(t)$ is an integrable
function
, defined almost everywhere in $T$, 
given by
    \begin{align}
        u_{\mathrm{eq}}(t) = \left(BK\right)^{-1} \frac{d}{dt} \sigma(t).
    \end{align}
\end{definition}

The next theorem provides the new condition for robust VSC design.

\begin{theorem}\label{thm:1}
    Consider the uncertain system~\eqref{eq:plant} and the sliding-mode controller~\eqref{eq:controller}. Given $\xi > 0$,
    if there exist diagonal matrices $W \in \mathbb{R}^{n \times n}$ and $X \in \mathbb{R}^{n \times n}$, a symmetric matrix $R \in \mathbb{R}^{n \times n}$, and a full matrix $Z \in \mathbb{R}^{m \times n}$, such that the following conditions hold:
    \begin{align}
        &W > 0, \quad R > 0, \label{eq:thm-1} \\
        &\begin{bmatrix}
            B_i Z + Z^\top B_i^\top + R & W - X + \xi Z^\top B_i^\top \\
            W - X + \xi B_i Z & - 2 \xi X
        \end{bmatrix} < 0,
        \forall i \in \mathbb{N}_{\leq{N}},
        \label{eq:thm-2}
    \end{align}
    then, the origin of the closed-loop system~\eqref{eq:closed-loop} with $K = Z X^{-1}$ is globally asymptotically stable. 
\end{theorem}
\begin{pf}
Consider that the conditions in~\eqref{eq:thm-1}--\eqref{eq:thm-2} hold. From~\eqref{eq:thm-2}, it follows that $X$ is a nonsingular matrix and there exists $X^{-1}$, since $X > 0$.
By multiplying the inequalities in~\eqref{eq:thm-2} by $\mathrm{diag}(X^{-1}, X^{-1})$ on the left and its transpose on the right, it follows that
\begin{align}\label{eq:proof-thm1-1}
    \begin{bmatrix}
            \Theta_i & P - X^{-1} + \xi K^\top B_i^\top X^{-1} \\
            P - X^{-1} + \xi X^{-1} B_i K & - 2 \xi X^{-1} 
        \end{bmatrix} < 0, 
\end{align}
for all $i \in \mathbb{N}_{\leq{N}}$, where $\Theta_i = X^{-1} B_i K + K^\top B_i^\top X^{-1} + Q$, $P = X^{-1} W X^{-1}$, $K = ZX^{-1}$, and $Q = X^{-1}RX^{-1}$. By multiplying~\eqref{eq:proof-thm1-1} by $[I \quad K^\top B_i^\top]$ on the left and its transpose on the right, it follows that
\begin{align}
    \label{eq:proof-thm1-2-1}
    P B_i K  + K^\top B_i^\top P + Q < 0, \quad i \in \mathbb{N}_{i \leq N}.
\end{align}
Since $B \in \mathrm{co}\{B_i\}_{i=1}^N$, if we multiply
\eqref{eq:proof-thm1-2-1} by $\alpha_i$ and sum all the inequalities from $1$ to $N$, we get 
\begin{align}
    \label{eq:proof-thm1-2}
    P B K  + K^\top B^\top P + Q < 0.
\end{align}

By multiplying~\eqref{eq:proof-thm1-2} with
$L(\sigma)$ on the left and the right, it follows that
\begin{align}\label{eq:proof-thm1-3}
   L(\sigma) P B K L(\sigma)  + L(\sigma) K^\top B^\top P L(\sigma) <  - L(\sigma)QL(\sigma) .
\end{align}
By multiplying~\eqref{eq:proof-thm1-3} by $x^\top$ on the left and its
transpose on the right, we obtain
\begin{align}\label{eq:proof-thm1-4}
    u_\mathrm{eq}^\top(t) PBK u_\mathrm{eq}(t) < - \frac{1}{2}u_\mathrm{eq}^\top(t) Q u_\mathrm{eq}(t),
\end{align}
since $u_\mathrm{eq}(t) = \mathrm{sgn}(\sigma(t)) =  L(\sigma(t)) x(t)$, where $u_\mathrm{eq}(t)$ is the extended control law given in Definition~\ref{def:extended_u}. Notice that~\eqref{eq:proof-thm1-2} implies that all eigenvalues of $BK$ lie
on the left half-plane, which ensures that $\mathrm{det}(BK) \neq 0$ and
$BK$ admits inverse. This guarantees that $u_{\mathrm{eq}}$ can be computed.

Since $V(\sigma)$ in~\eqref{eq:proof-thm1-6} is a Lipschitz function and $\sigma(t)$ is
absolutely continuous, then $V(\sigma(t))$ is
also absolutely continuous. It ensures the existence 
of the time-derivative of $V(\sigma(t))$ along 
a Filippov solution of~$\sigma(t)$. 
As argued by~\cite{hsu2000matrix}, if $V(\sigma)$ were continuously differentiable, then
\begin{align}\label{eq:proof-thm1-7}
    \frac{dV}{dt} = \left(\frac{\partial V}{\partial \sigma}(t)\right)^\top \dot{\sigma}(t) = P B K u_{\mathrm{eq}}(t),
\end{align}
However, the computation in~\eqref{eq:proof-thm1-7} can not be performed because $\frac{\partial V}{\partial \sigma}(t) = p_i \mathrm{sgn}(\sigma_i(t))$ is not defined whenever $\sigma_i(t) = 0$, for some $i \in \mathbb{N}_{\leq N}$.
In order to properly compute $\dot{V} = dV/dt$, we need to show that
\begin{align}\label{eq:proof-thm1-8}
    \frac{\partial V}{\partial \sigma}(t) = P u_\mathrm{eq}(t)
\end{align}
holds almost everywhere. For that purpose, consider that $\sigma_i \neq 0$, 
for all $i \in \mathbb{N}_{\leq N}$. Then, 
\begin{align}
    \left(\frac{\partial V}{\partial \sigma_i}(t)\right) = p_i \mathrm{sgn}(\sigma_i(t)) = p_i u_{\mathrm{eq}_i}(t).
\end{align}
Now, consider some open interval $\mathcal{T}$ of $t$. 
Then, there exists $i \in \mathbb{N}_{\leq N}$ such that
$\sigma_i(t) \equiv 0$, which implies that $\dot{\sigma}_i(t) \equiv 0$ in $\mathcal{T}$ and
\begin{align}
    \left(\frac{\partial V}{\partial \sigma_i}(t)\right) \dot{\sigma}_i(t) = p_i u_{\mathrm{eq}_i}(t) \dot{\sigma_i}(t) \equiv 0
\end{align}
holds in $\mathcal{T}$. Since all the other possibilities include the case of some $\sigma_i = 0$ only in sets of measure zero, it is possible to conclude that~\eqref{eq:proof-thm1-8} holds almost everywhere.
Since 
\begin{align}\label{eq:proof-thm1-9}
    \dot{V}(t) = \sum_{i=1}^n \left(\frac{\partial V}{\partial \sigma_i}(t)\right) \dot{\sigma_i}(t),
\end{align}
we can conclude from~\eqref{eq:closed-loop},~\eqref{eq:proof-thm1-4}, and~\eqref{eq:proof-thm1-9} that
\begin{align}\label{eq:proof-thm1-10}
    \dot{V}(t) = u_\mathrm{eq}^\top(t) PBK u_\mathrm{eq}(t) < - \frac{1}{2}u_\mathrm{eq}^\top(t) Q u_\mathrm{eq}(t) < 0,
\end{align}
provided that $Q$ is positive definite and $u_\mathrm{eq}(t)$ is null only when $\sigma(t) \equiv 0$. 
It implies that $\dot{V}(\sigma(t))$ converges to zero in finite time
and the origin is globally attractive. In fact, $\|\mathrm{sgn}(\sigma)\|_2 \geq 1$, $\forall \sigma \neq 0$. Since $P > 0$, it ensures that $V(\sigma)$ is positive definite.
Although $V(\sigma)$ is not differentiable, it is continuous and lower bounded by a class $\mathcal{K}$ function. Those properties allow to conclude the global stability since $dV(\sigma(t))/dt < 0$, 
$\forall \sigma \neq 0$, almost everywhere.
To show the finite-time convergence, notice that $\|\mathrm{sgn}(\sigma)\|_2 \geq 1$, $\forall \sigma \neq 0$.
Then, we can obtain that
$V(t) \leq V(0) -\frac{1}{2}\lambda_{\min}(Q) t$.
If all states $\sigma_i$, $i \in \mathbb{N}_{\leq n}$, achieve the sliding surface at $t=T_{\mathrm{vsc}}$, one has that $V(t_\mathrm{vsc})=0$ and
\begin{align}\label{eq:reaching_time}
    t_{\mathrm{vsc}} \leq 2V_0/\lambda_{\min}(Q),
\end{align}
where $V_0 = V(\sigma(0)) = \sum_{i=1}^n p_i |\sigma_i(0)|$.
This concludes the proof.
\end{pf}

\subsection{Optimization issues related to the reaching time of VSC}
\label{sec:optimization_VSC}
In this section, we provide optimization procedures to estimate the upper bound of the
reaching time and evaluating the influence of the initial condition and the convergence
associated with the decay rate of the Lyapunov function~\eqref{eq:proof-thm1-6}.

From~\eqref{eq:reaching_time}, it is possible to notice that
the upper bound for the reaching time depends on the
initial condition $V_0$, 
and the smallest eigenvalue of the positive definite matrix $Q$, which is related
to the decay rate of ${V}(t)$. 
For a given initial condition $\sigma(0)$ (associated to $V_0$), 
the reaching time can be minimized by maximizing the
smallest eigenvalue of $Q$. This objective can be achieved by incorporating the following
constraint:
\begin{align}\label{eq:maximize_Q}
\begin{bmatrix}
    R & X \\
    X & \rho_\mathrm{vsc} I
\end{bmatrix} \geq 0.
\end{align}
From~\eqref{eq:maximize_Q}, it follows from Schur complement that
$R - \rho_\mathrm{vsc}^{-1} X^2 \geq 0$. By multiplying both sides by $X^{-1}$,
we have that $Q \geq \rho^{-1} I$, since $Q = X^{-1} R X^{-1}$.
Thus, by minimizing $\rho$, the eigenvalues of $Q$ are maximized, thus
reducing the reaching time $t_{\mathrm{vsc}}$.

To evaluate the influence of the initial condition on the
reaching time $t_{\mathrm{vsc}}$, we consider the following constraint:
\begin{align}\label{eq:maximize_V0}
    \begin{bmatrix}
        \varphi_\mathrm{vsc} I & I \\ I & 2X-W
    \end{bmatrix} \geq 0.
\end{align}
Inequality~\eqref{eq:maximize_V0} implies from Schur complement that $\varphi I \geq (2X-W)^{-1}$.
However, since $(P^{-1}-X)P(P^{-1}-X) \geq 0$, provided that $P > 0$, then $P^{-1} \geq 2X-W$, since $W = XPX$.
Thus, $P \leq (2X-W)^{-1}$, which implies that $P \leq \varphi_\mathrm{vsc} I$ or still $V(x) \leq \varphi_\mathrm{vsc} x^\top x$.
Hence, it is possible to conclude that $\mathcal{B}_\mathrm{vsc} \subset \Omega_{\mathrm{vsc}}$,
where $\mathcal{B}_\mathrm{vsc} = \{x \in \mathbb{R}^n : x^\top x \leq \varphi_\mathrm{vsc}^{-1}\}$ and 
\begin{align}\label{eq:Omega-set}
    \Omega_{\mathrm{vsc}} = \{\sigma \in \mathbb{R}^n : V(\sigma) \leq 1\}.
\end{align}
Thus, if $\varphi_\mathrm{vsc}$ is minimized, the set $\Omega_{\mathrm{vsc}}$ is enlarged. Therefore, if $\sigma(0)$ is taken inside of $\Omega_\mathrm{vsc}$, the reaching time is bounded by 
\begin{align}
t_{\mathrm{vsc}} \leq 2V_0/\lambda_{\min}(Q) \leq 2{\rho}_\mathrm{vsc} = T_{\mathrm{vsc}}. 
\end{align}
To evaluate the trade-off between the size of the set of initial conditions $\Omega_{\mathrm{vsc}}$ and the convergence rate associated with $\lambda_{\min}(Q)$, we formulate a convex optimization problem
to obtain the minimum reaching time associated with a given initial condition set. 
The optimization problem is stated in the sequel, for a given $\varphi_\mathrm{vsc} >0$:
\begin{flalign}
&\min\limits_{R, W, X, Z} \qquad \rho_\mathrm{vsc}  \label{eq:optmization_problem-1}\\
&\text{subject to}~\text{LMIs in}~\eqref{eq:thm-1}, \eqref{eq:thm-2}, \eqref{eq:maximize_Q}, \eqref{eq:maximize_V0}. \nonumber 
\end{flalign}




\section{Robust Unit Vector Control Design}
\label{sec:UVC}
Consider now the following unit vector control law:
\begin{align}\label{eq:control-uvc}
    u = K \frac{\sigma}{\|\sigma\|},
\end{align}
which leads to the following closed-loop system obtained by substituting~\eqref{eq:control-uvc} in~\eqref{eq:plant}:
\begin{align}\label{eq:closed-loop-uvc}
    \dot{\sigma} = BK \frac{\sigma}{\|\sigma\|}.
\end{align}
Let $z = R(\sigma) \sigma$, where $R(\sigma) = 1/\sqrt{\|\sigma\|}$. In $z$-coordinates, the closed-loop system~\eqref{eq:closed-loop-uvc}
can be equivalently rewritten as
\begin{align}
    \dot{z} = -\frac{1}{2}R(\sigma) \Pi_\sigma BK R(\sigma) z + R(\sigma) BK R(\sigma) z,
\end{align}
where $\Pi_\sigma = \sigma \sigma^\top/\|\sigma\|^2$ is a projection matrix which satisfies the
following properties: $\Pi_\sigma = \Pi_\sigma^\top$, $\Pi_\sigma^2 = \Pi_\sigma$, and $\|\Pi_\sigma\| = 1$, $\forall \sigma \in \mathbb{R}^n$. 

The problem addressed in this section is stated as follows.
\begin{problem}\label{problem:2}
    Consider the uncertain system~\eqref{eq:plant} and the UVC law~\eqref{eq:control-uvc}.
    Design a robust control gain $K$ such that the origin of the closed-loop system~\eqref{eq:closed-loop} is globally finite-time stable. 
\end{problem}

\subsection{The proposed robust UVC design condition}

The robust design condition for the UVC is derived in this section using the following Lyapunov function candidate:
\begin{align}\label{eq:Lyap-uvc-1}
    U(\sigma) = \frac{1}{\|\sigma\|} \sigma^\top P \sigma,
\end{align}
where $P = P^\top > 0$. This function can be rewritten as the following standard quadratic function using the $z$-coordinates:
\begin{align}\label{eq:Lyap-uvc-2}
    U(z) = z^\top P z.
\end{align}

The next theorem provides the new condition for robust UVC design.
\begin{theorem}\label{thm:2}
 Consider the uncertain system~\eqref{eq:plant} and the sliding-mode controller~\eqref{eq:controller}. Given $\mu > 0$,
    if there exist symmetric matrices $X \in \mathbb{R}^{n \times n}$ and $R \in \mathbb{R}^{n \times n}$, and a full matrix $Z \in \mathbb{R}^{m \times n}$, such that the following conditions hold:
    \begin{align}
        &X > 0, \quad R > 0, \label{eq:thm-1-2} \\
        &\begin{bmatrix}
             B_i Z + Z^\top B_i^\top + \dfrac{\mu}{4} I + R & Z^\top B_i^\top \\
            B_i Z & -\mu I
        \end{bmatrix} < 0,
        \forall i \in \mathbb{N}_{\leq{N}},
        \label{eq:thm-2-2}
    \end{align}
    then, the origin of the closed-loop system~\eqref{eq:closed-loop-uvc} with $K = Z X^{-1}$ is globally asymptotically stable.     
\end{theorem}
\begin{pf}
    Consider that the conditions in~\eqref{eq:thm-1-2}--\eqref{eq:thm-2-2} hold. From~\eqref{eq:thm-1-2}, it follows that $X$ is a nonsingular matrix and there exists $X^{-1}$, since $X > 0$.
By multiplying the inequalities in~\eqref{eq:thm-2-2} by $\mathrm{diag}(X^{-1}, I)$ on the left and its transpose on the right, it follows that
\begin{align}\label{eq:proof-thm2-1}
    \begin{bmatrix}
            PB_iK + K^\top B_i^\top P + Q + \dfrac{\mu}{4} P^2 & K^\top B_i^\top \\
            B_i K & -\mu I
        \end{bmatrix} < 0, 
\end{align}
for all $i \in \mathbb{N}_{\leq{N}}$, where $P = X^{-1}$, $K = ZX^{-1}$, and $Q = X^{-1}RX^{-1}$. 
Since $B \in \mathrm{co}\{B_i\}_{i=1}^N$, if we multiply
\eqref{eq:proof-thm2-1} by $\alpha_i$ and sum all the inequalities from $1$ to $N$, and then apply Schur complement, we can obtain
\begin{align}
    \label{eq:proof-thm2-2}
    \frac{1}{\mu} K^\top B^\top B K + \frac{\mu}{4} P^2 + P B K + K^\top B^\top P + Q < 0.
\end{align}
Provided that 
\begin{align}    \label{eq:proof-thm2-3}
    - \frac{1}{2}K^\top B^\top \Pi_\sigma P - \frac{1}{2}P \Pi_\sigma BK \leq \frac{1}{\mu} K^\top B^\top B K + \frac{\mu}{4} P^2 
\end{align}
since
\begin{align*}
    \left(\frac{1}{\sqrt{\mu}}BK + \frac{\sqrt{\mu}}{2}\Pi_\sigma  P \right)^\top \left(\frac{1}{\sqrt{\mu}}BK + \frac{\sqrt{\mu}}{2}\Pi_\sigma P \right) \geq 0
\end{align*}
and $\|\Pi_\sigma\|=1$, then it follows from~\eqref{eq:proof-thm2-2} and \eqref{eq:proof-thm2-3} that
\begin{align}\label{eq:proof-thm2-4}
    - \frac{1}{2}P \Pi_\sigma BK - \frac{1}{2}K^\top B^\top \Pi_\sigma P \!+\! S BK \!+\! K^\top B^\top P + Q < 0.
\end{align}

By multiplying~\eqref{eq:proof-thm2-4} with
$z^\top R(\sigma)$ on the left and its transpose on the right, it follows that
\begin{align}\label{eq:proof-thm2-5}
   \dot{U}(z) < - z^\top R(\sigma)QR(\sigma)z < 0,
\end{align}
with $U(z)$ defined in~\eqref{eq:Lyap-uvc-2}. By following similar arguments of the proof of Theorem~\ref{thm:1}, it is possible to ensure that the origin is globally attractive. 
To show the finite-time convergence, notice that 
$z^\top R(\sigma) Q R(\sigma) z \geq \lambda_{\min}(Q) {\|z\|^2}/{\|\sigma\|} =  \lambda_{\min}(Q)$,
hence, it is possible to obtain from \eqref{eq:proof-thm2-5} that the reaching time is upper-bounded by
\begin{align}
    t_{\mathrm{uvc}} \leq U_0/\lambda_{\min}(Q),
\end{align}
where $U_0 = U(\sigma(0)) = \sigma^\top(0) P \sigma(0) / \|\sigma(0)\|$, $\forall \sigma(0) \neq 0$.
This concludes the proof.
\end{pf}

\subsection{Optimization issues related to the reaching time of UVC}

In this section, we employ the same reasoning as in Section~\ref{sec:optimization_VSC} to formulate convex optimization problems for designing the robust UVC. 
For a given initial condition $\sigma(0)$ (associated to $U_0$), 
the reaching time can also be minimized by maximizing the
smallest eigenvalue of $Q$ in Theorem~\ref{thm:2}. This objective can be achieved by~\eqref{eq:maximize_Q}. However, notice that $X = P^{-1}$ in Theorem~\ref{thm:2}.
To enlarge the estimated set of initial conditions, we consider the following constraint:
\begin{align}\label{eq:maximize_U0}
    \begin{bmatrix}
        \varphi_\mathrm{uvc} I & I \\ I & X
    \end{bmatrix} \geq 0.
\end{align}
Inequality~\eqref{eq:maximize_U0} implies from Schur complement that $\varphi I \geq X^{-1}$.
Thus, $P \leq \varphi I$ or still $U(z) \leq \varphi z^\top z$.
Hence, it is possible to conclude that $\mathcal{B}_\mathrm{uvc}\subset \Omega_{\mathrm{uvc}}$,
where $\mathcal{B}_\mathrm{uvc} = \{z \in \mathbb{R}^n : z^\top z \leq \varphi^{-1}\}$ and 
\begin{align}\label{eq:Omega-set-2}
    \Omega_{\mathrm{uvc}} = \{\sigma \in \mathbb{R}^n : U(\sigma) \leq 1\}.
\end{align}
Thus, if $\varphi$ is minimized, the set $\Omega_{\mathrm{uvc}}$ is enlarged. Therefore, if $\sigma(0)$ is taken inside of $\Omega$, the reaching time is bounded by 
\begin{align}
t_{\mathrm{uvc}} \leq V_0/\lambda_{\min}(Q) \leq {\rho}_{\mathrm{uvc}} = T_{\mathrm{uvc}}. 
\end{align}

The optimization problem for minimizing the reaching time estimate for a given set of initial conditions associated with $\varphi_\mathrm{uvc} >0$ is stated in the sequel:
\begin{flalign}
&\min\limits_{R, X, Z} \qquad \rho_{\mathrm{uvc}} \label{eq:optmization_problem-3}\\
&\text{subject to}~\text{and LMIs in}~\eqref{eq:thm-1-2}, \eqref{eq:thm-2-2}, \eqref{eq:maximize_Q}, \eqref{eq:maximize_U0}. \nonumber 
\end{flalign}


\section{Numerical Results}
\label{sec:results}

This section presents two numerical examples to illustrate
the effectiveness of the proposed control design conditions.
The first example considers the second-order dynamical model of a 
robotics visual servo problem. The second example 
considers an over-actuated underwater Remotely Operated Vehicle (ROV) model.

\subsection{Example~1: Robotics visual servo system}

Consider a planar kinematic manipulator with an end-effector 
image position coordinates $\sigma = [p_x,p_y]^\top \in \mathbb{R}^2$ 
given by an uncalibrated fixed camera with an optical orthogonal axis concerning the
robot workspace plane~\cite{oliveira2014overcoming,oliveira2010sliding}. The uncertain dynamics are described by
\begin{align}
    B(\phi) = 
    \begin{bmatrix}
        \cos{(\phi)} & \sin{(\phi)} \\
        -\sin{(\phi)} & \cos{(\phi)}
    \end{bmatrix},
\end{align}
which is a matrix that depends on the uncertain rotation angle $\phi$
due to the uncalibrated camera. Let $\bar{\phi}$ be a given nominal angle,
the uncertainty can be modeled as the variation $\Delta \phi = \phi - \bar{\phi}$,
such that $|\Delta \phi| \leq \bar{\Delta}$ and $B(\phi) = B(\Delta \phi) B(\bar{\phi})$. With this uncertainty description, it is possible to obtain
$N=4$ vertices for the polytopic representation of the uncertain matrix $B(\phi)$ associated with the variations of the following parameters:
\begin{align}
    \begin{bmatrix}
        \cos{(\Delta \phi)} \\
        \sin{(\Delta \phi)}
    \end{bmatrix}
    \in
    \mathrm{co}\left\lbrace 
    \begin{bmatrix}
        \cos{(\Delta \phi)} \\
        \sin{(\Delta \phi)}
    \end{bmatrix},
    \begin{bmatrix}
        \cos{(\Delta \phi)} \\
        \sin{(\Delta \phi)}
    \end{bmatrix}, \right. \nonumber \\
    \left.
    \begin{bmatrix}
        \cos{(\Delta \phi)} \\
        \sin{(\Delta \phi)}
    \end{bmatrix},
    \begin{bmatrix}
        \cos{(\Delta \phi)} \\
        \sin{(\Delta \phi)}
    \end{bmatrix}
    \right\rbrace,
\end{align}
valid for all $|\Delta \phi| \leq \bar{\Delta}$, provided that $0 \leq \bar{\Delta} \leq \pi/2~\mathrm{rad}$. 
For the conducted experiments, we assume that $\bar{\phi} = \pi/6$ and $\bar{\Delta} = \pi/4$.

\begin{figure}[!ht]
    \centering
    \includegraphics[width=0.8\columnwidth]{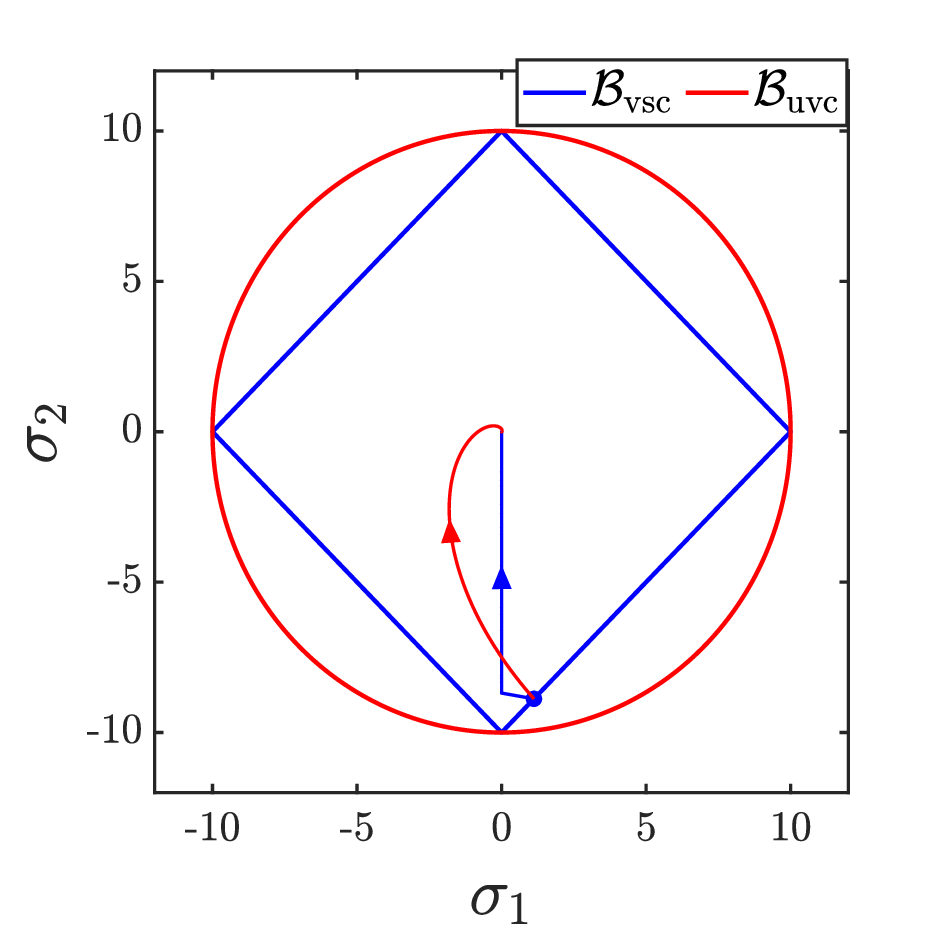}
    \caption{Regions $\mathcal{B}_{\mathrm{vsc}}$ and $\mathcal{B}_{\mathrm{uvc}}$ and closed-loop trajectories with the VSC (in blue) and with the UVC (in red) -- Example~1.}
    \label{fig:pplane_ex1}
\end{figure}

To fairly compare the design conditions for the VSC and UVC, we 
fix $\rho_{\mathrm{vsc}} = 0.25$, $\rho_{\mathrm{uvc}} = 0.5$, and $\varphi_{\mathrm{vsc}} = \varphi_{\mathrm{uvc}} = 0.1$.
Notice that both $T_{\mathrm{vsc}} = 2{\rho}_{\mathrm{vsc}} = 0.5$~s and $T_{\mathrm{uvc}} = \rho_{\mathrm{uvc}} = 0.5$~s.
Based on these values, we can solve feasibility conditions to design the robust VSC and UVC using optimization problems based on Theorem~\ref{thm:1} and Theorem~\ref{thm:2}, respectively.
Feasible solutions were obtained with $\xi = 0.001$ in Theorem~\ref{thm:1} and $\mu = 1000$ in Theorem~\ref{thm:2}. The designed (non-diagonal) control gains are:
\begin{align*}
    K_{\mathrm{vsc}} &= 
    \begin{bmatrix}
        -33.2438 & 19.1933 \\
        -19.1933 & -33.2438
    \end{bmatrix},\\
    K_{\mathrm{uvc}} &= 
    \begin{bmatrix}
        -52.8970 &  30.5401 \\
        -30.5401 & -52.8970
    \end{bmatrix}.
\end{align*}
The regions $\mathcal{B}_{\mathrm{vsc}}$ and $\mathcal{B}_{\mathrm{uvc}}$ are depicted in Fig.~\ref{fig:pplane_ex1} in the $\sigma$-coordinates. It is possible to notice that, due to the induced norms, the regions have different shapes and  $\mathcal{B}_{\mathrm{uvc}}$ contains $\mathcal{B}_{\mathrm{vsc}}$.

\begin{figure}[!ht]
    \centering
    \includegraphics[width=\columnwidth]{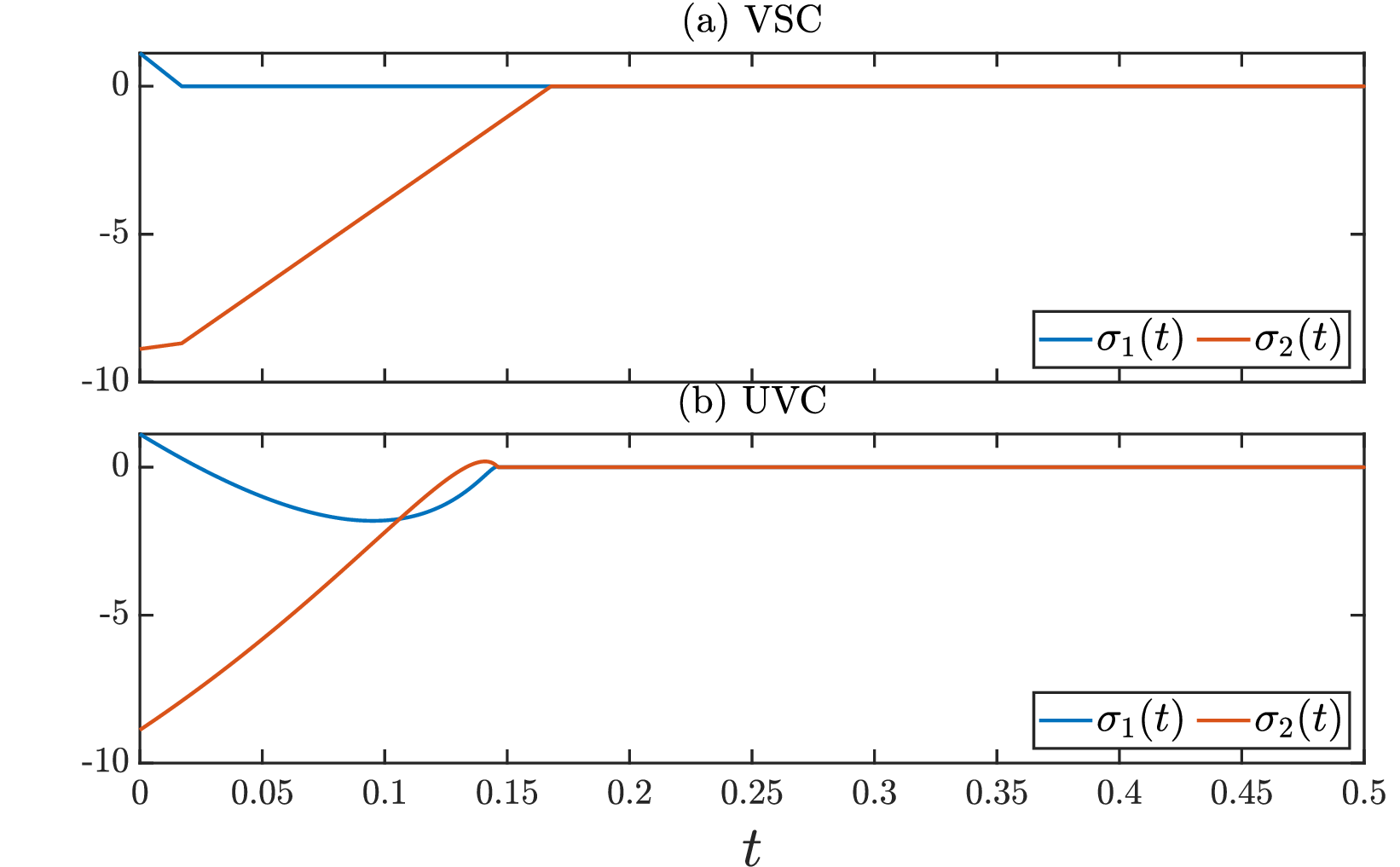}
    \caption{States of the closed-loop robotics visual servo system with the robust (a) VSC and (b) UVC -- Example~1.}
    \label{fig:states_ex1}
\end{figure}

The states of the closed-loop system with the VSC and the UVC are shown in Fig.~\ref{fig:states_ex1}, for $B = B_3$. As observed in Figs.~\ref{fig:pplane_ex1} and~\ref{fig:states_ex1}, the states of the closed-loop system with the VSC can converge to the sliding surface independently, while in the case of the UVC, the sliding surface is achieved at the origin~\cite{hsu2002multivariable}. However, in both cases, the states reach the sliding surface in a time smaller than $T_\mathrm{vsc} = T_\mathrm{uvc} = 0.5$~s.

\subsection{Example~2: Underwater ROV system}

The underwater ROV system, borrowed from~\cite{geromel2024multivariable}, is described with the state $\sigma=[v_x \; v_y \; \omega_z]^\top \in \mathbb{R}^3$, where $v_x$ and $v_y$ are velocities related to the body frame and $\omega_z$ is the angular velocity with respect to the $z$-axis. The system has four inputs resulting from propellers responsible for the displacement of the body. According to~\cite{geromel2024multivariable}, a simplified uncertain model of this system can be obtained by taking $B(g)=M^{-1}\Psi\Pi(g)$, with  $M=\mathrm{diag}(m_0, m_o, I_z)$, and
\begin{align}
    \Psi =
    \begin{bmatrix}
        \psi_1 & \psi_1 & \psi_1 &  \psi_1 \\
        \psi_1 & -\psi_1 & -\psi_1 & \psi_1\\
        -\psi_2 & \psi_2 & -\psi_2 & \psi_2
    \end{bmatrix},
\end{align}
where $m_0 = 290~\mathrm{kg}$ is the ROV mass, $I_z = 290~\mathrm{kg m}^2$ is the
moment of inertia, $\psi_1 = \sqrt{2}/2$, and $\psi_2 = 0.35~\mathrm{m}$. 
The input matrix with uncertain coefficients is $\Pi(g) = \mathrm{diag}(g_1,1,g_3,1)$,
where $g_1, g_3 \in [1/2,1]$ are uncertain gains in the actuator channels. 
Since the system has two uncertain parameters, it leads to $N=4$ vertices $B_i \in \mathbb{R}^{3 \times 4}$ in the polytopic description of the matrix $B$ in~\eqref{eq:plant}.

\begin{figure}[!ht]
    \centering
    \includegraphics[width=\columnwidth]{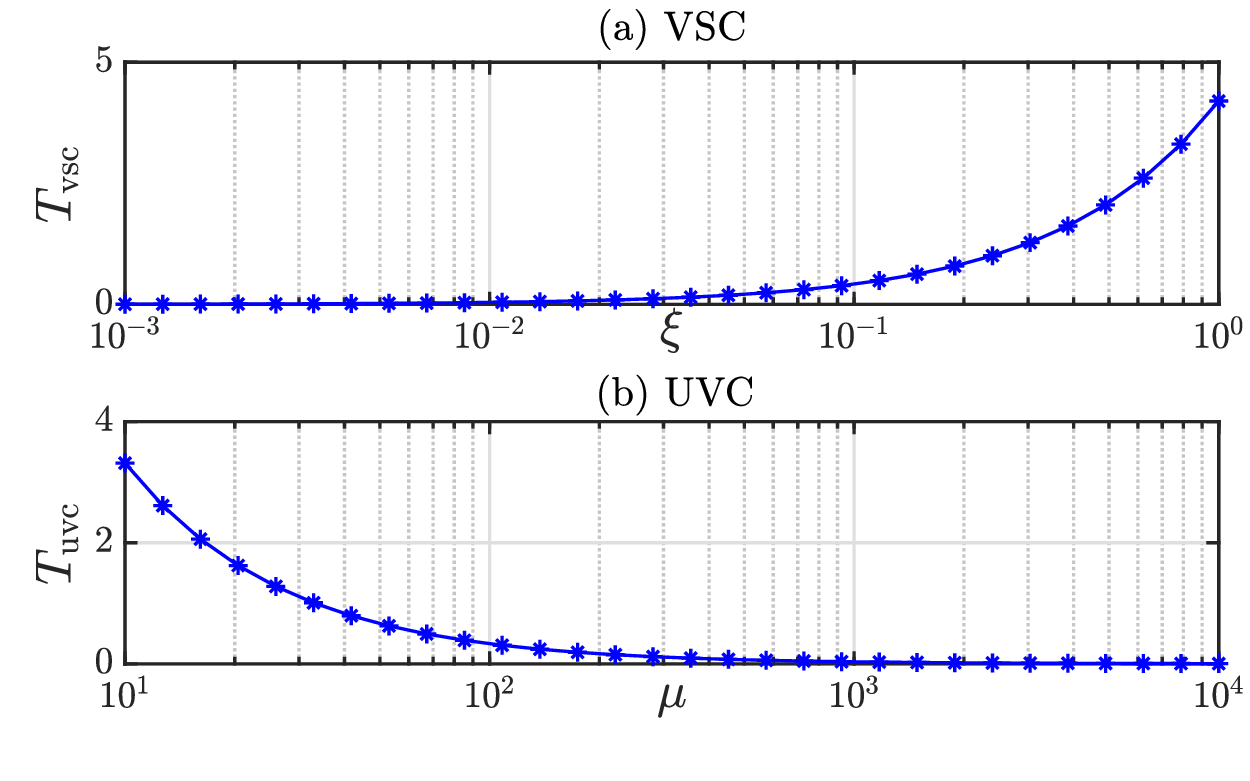}
    \caption{Upper-bound on the reaching times for the (a) VSC and (b) UVC -- Example~2.}
    \label{fig:Tr_ex2}
\end{figure}

Consider a given initial condition~$x(0) = [1 \; 1 \; \pi/4]^\top$. 
To evaluate the influence of the scalar parameters $\xi$ and $\mu$ over the reaching time, 
we solve optimization problems~\eqref{eq:optmization_problem-1} 
and~\eqref{eq:optmization_problem-3}, respectively, 
for different values of parameters $\xi$ and $\mu$ and $\varphi_{\mathrm{vsc}} = \varphi_{\mathrm{uvc}} = 0.4$, which ensures the inclusion of the given initial condition within the guaranteed sets of initial conditions. The results are depicted in Fig.~\ref{fig:Tr_ex2}. 
Notice that, in this example, smaller reaching times can be obtained with smaller values of $\xi$ in optimization problem~\eqref{eq:optmization_problem-1} and larger values of $\mu$
in optimization problem~\eqref{eq:optmization_problem-3}.

For $\xi = 0.2395$ and $\mu = 32.9034$, we obtain $T_{\mathrm{vsc}} = 1.0057$~s,
$T_{\mathrm{uvc}} = 1.0077$~s, and
\begin{align*}
    K_{\mathrm{vsc}} &= \begin{bmatrix}
        -677.59 &  2.42 &  291.63 \\
         -469.25 & 1238.52 & -123.22 \\
         -677.58 & -2.42 &  291.64 \\
         -469.25 & -1238.52 & -123.22
    \end{bmatrix}, \\
    K_{\mathrm{uvc}} &= 
    \begin{bmatrix}
        -874.02 & -0.39 &  227.06 \\
         -598.96 &   1299.41 & -139.86 \\
         -874.03 &  0.39 &  226.97\\
         -598.92 &  -1299.43 &  -139.6
    \end{bmatrix}.
\end{align*}

\begin{figure}[!ht]
    \centering
    \includegraphics[width=\columnwidth]{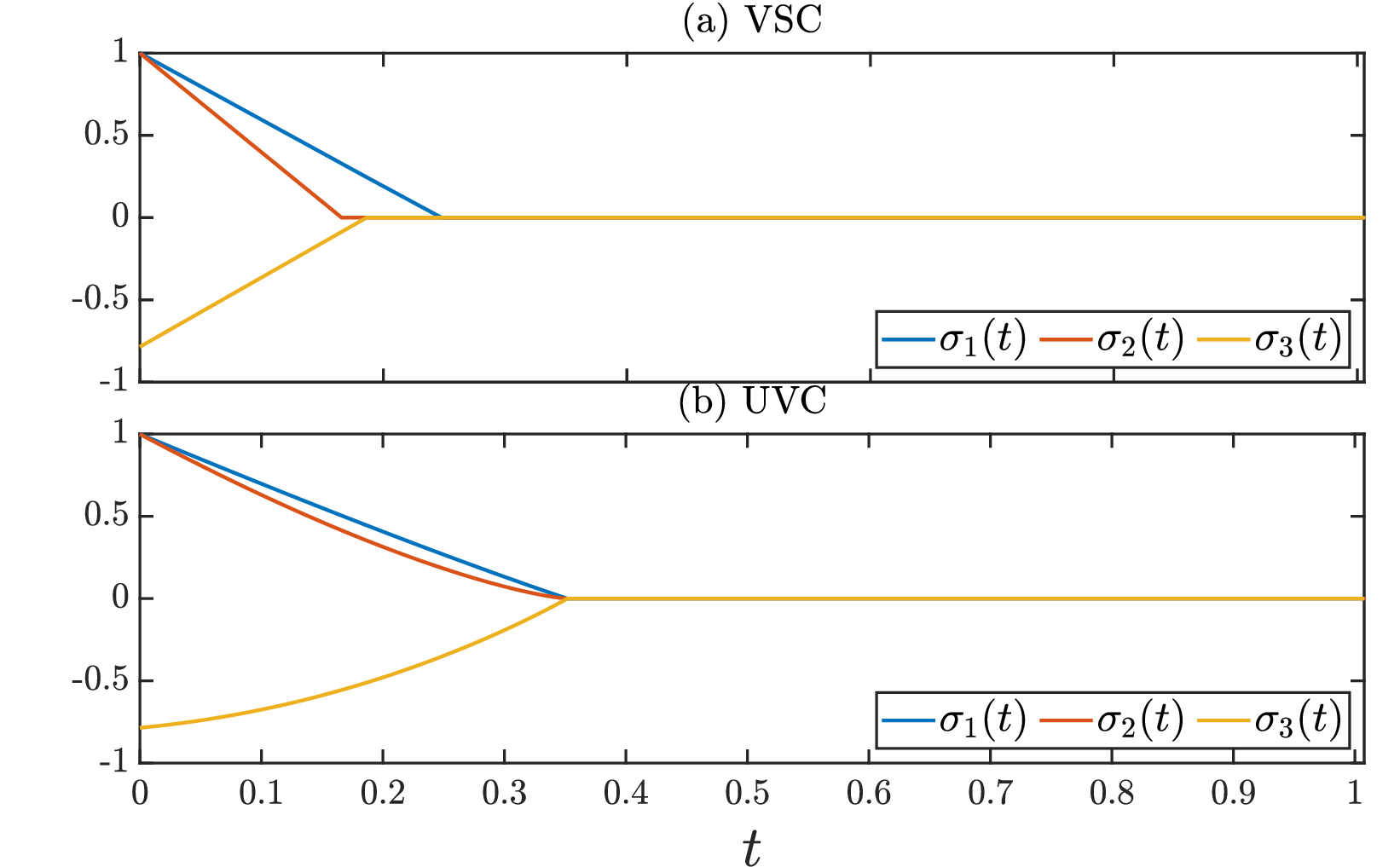}
    \caption{States of the closed-loop over-actuated ROV system with the robust (a) VSC and (b) UVC -- Example~2.}
    \label{fig:states_ex2}
\end{figure}

The trajectories of the closed-loop system with the VSC and the UVC are
depicted in Fig.~\ref{fig:states_ex2} for $B = B_1$.


\section{Conclusion}
\label{sec:conclusion}

This paper has addressed the multivariable sliding mode control design for a class of polytopic uncertain systems. Based on suitable representations of the closed-loop system, we developed LMI-based conditions to design robust VSC and UVC laws to ensure the global stability of the closed-loop system. Moreover, by noticing that the reaching time depends on the initial condition and the decay rate, we provided optimization problems to design the controllers by taking into account the trade-off between the initial condition and the decay rate. 
This optimization procedure allows the minimization of the reaching time related to a given set of initial conditions. The results indicated the effectiveness of the proposed approach in designing the globally stabilizing robust sliding mode controllers. An interesting aspect is that the state can converge to the sliding surface independently with the VSC, different from UVC in which the sliding mode occurs only at the origin. As future research, we intend to investigate conditions for systems in the presence of matched and unmatched disturbances~\cite{ovalle2024sliding}.

\begin{ack}                               
This work was supported by the Brazilian agencies CNPq (Grant numbers: 407885/2023-4 and 309008/2022-0), CAPES and FAPERJ.  
\end{ack}

\bibliographystyle{abbrv}
\bibliography{references}           



\end{document}